\documentclass[11pt]{amsart}
\usepackage{amssymb,amsbsy,latexsym,amsmath,bbm,epsfig,psfrag,amsthm,mathrsfs}
\usepackage[swedish,english]{babel}
\usepackage{graphicx,tikz,esint}
\usepackage[all,cmtip]{xy}
\usepackage[T1]{fontenc}
\usepackage[latin1]{inputenc}
\usepackage{amsfonts}
\usepackage{amsmath}
\usepackage{amssymb}
\usepackage{epsf}
\usepackage{graphics}
\usepackage{graphicx}
\usepackage{latexsym}
\usepackage{psfrag}
\usepackage{relsize}
\usepackage{verbatim}
\usepackage{hyperref}
\usepackage{pgfplots}
\usepackage{float}
\usetikzlibrary{positioning}
\usetikzlibrary{arrows}
\usetikzlibrary{decorations.pathreplacing}
\usepackage{slashed}

\usepackage{url}

\theoremstyle{plain}

\newtheorem{Lem}{Lemma}
\newtheorem{Prop}{Proposition}
\newtheorem{Thm}{Theorem}
\numberwithin{Con}{section}

\theoremstyle{definition}
\newtheorem{Def}{Definition}
\numberwithin{Qu}{section}

\numberwithin{hyp}{section}

\numberwithin{conj}{section}
\newtheorem{ex}{Example}
\numberwithin{ex}{section}

\theoremstyle{remark}
\newtheorem{rem}{\bf{Remark}}
\numberwithin{rem}{section}

\numberwithin{equation}{section}

\DeclareMathOperator*{\Div}{div}
\DeclareMathOperator*{\Ker}{ker}

\DeclareMathOperator*{\im}{im}

\newcommand{\dv}{\partial}

\newcommand{\Om}{\Omega}
\newcommand{\dbar}{{\overline \partial}}

\newcommand{\eps}{\varepsilon}

\newcommand{\R}{{\mathbb R}}
\newcommand{\C}{{\mathbb C}}

\newcommand{\Z}{{\mathbb Z}}

\newcommand{\N}{{\mathbb N}}

\newcommand{\Di}{\mathbb{D}}
\newcommand{\De}{\mathscr{D}}

\newcommand{\LL}{\mathcal{L}}
\newcommand{\Bb}{\mathbb{B}}
\newcommand{\A}{\mathscr{A}}
\newcommand{\Aa}{\mathbb{A}}
\newcommand{\Pop}{\mathscr{P}}
\newcommand{\Pp}{\mathbb{P}}

\newcommand{\Qq}{\mathbb{Q}}

\newcommand{\B}{\mathscr{B}}

\newcommand{\ri}{\,\,\lrcorner\,\,}

\newcommand{\zbar}{{\overline z}}
\newcommand{\wedg}{\mathbin{\scriptstyle{\wedge}}}
\newcommand{\lctr}{\mathbin{\lrcorner}}
\mathchardef\semic="303A
\newcommand{\sett}[2]{ \{ #1 \, \semic \, #2 \} }

\usepackage{color}

\definecolor{gr}{rgb}   {0.,   0.8,   0. }
\definecolor{bl}{rgb}   {0.,   0.5,   1. }
\definecolor{mg}{rgb}   {0.7,  0.,    0.7}

\begin{document}
%\begin{flushright}

%December 8, 2008 (first part) 
%January 9, 2009 (second part) 
%\end{flushright}
\vspace{.4cm}

\title[]
{\bf \sffamily Coerciveness and Morrey inequalities for elliptic operators with natural boundary conditions via Weitzenb\"ock identities}
\author[]{Erik Duse and Andreas Ros\'en}

\address{Erik Duse \\ Department of mathematics, KTH \\SE-100 44 Stockholm, Sweden} \email{duse@kth.se}
\address{Andreas Ros\'en\\Mathematical Sciences, Chalmers University of Technology and University of Gothenburg\\
SE-412 96 G{\"o}teborg, Sweden}
\email{andreas.rosen@chalmers.se}

%\keywords{}
%\subjclass[2010]{42B35, 42B20, 42B37}

%\date{}

\dedicatory{}

\thanks{ED was supported by the Swedish Research Council (VR), grant no. 2019-04152 and the Knut and Alice Wallenberg Foundation grant KAW 2015.0270.}

\thanks{AR was supported by Grant 2022-03996 from the Swedish research council, VR}

\maketitle

\begin{abstract}
We prove a Weitzenb\"ock identity for general pairs of constant coefficient homogeneous first order partial differential operators, and deduce from it
sufficient algebraic conditions for coerciveness and Morrey estimates under 
the natural 1/2 boundary conditions.
Our proof of the $W^{1,2}$ elliptic estimate relies on the Aronszajn-Ne{\u c}as-Smith
coercive estimate.
For generalized strongly pseudoconvex domains, we improve the Morrey estimate
to a weighted $W^{1,2}$ square function estimate, using a generalized
Cauchy--Pompieu reproducing formula and the $T1$ theorem for singular integrals. 
We use Van Schaftingen's notion of cocanceling to study the generalized Levi
forms appearing.
\end{abstract}

\section{\sffamily  Introduction}

In this work we will prove a Weitzenb\"ock identity for pairs of constant coefficient homogeneous first order partial differential operators on $C^2$ domains in $\R^n$, for functions that satisfy a natural local boundary condition. 
We denote by $E$, $F$ and $G$ real or complex inner product spaces. 
In the complex case, we shall denote by $\langle\cdot,\cdot\rangle$ the 
real part of the complex inner product.  
Given differential operators $\A:C^\infty(\R^n,F)\to C^\infty(\R^n,G)$ and $\B: C^\infty(\R^n, F)\to C^\infty(\R^n, E)$, then in euclidean coordinates of $\R^n$ the differential operators $\A$ and $\B$ are represented in the form
\begin{align*}
\A=\sum_{j=1}^nA_j\dv_j,\quad \B=\sum_{j=1}^nB_j\dv_j,
\end{align*}
where $A_j\in \LL(F,G)$ and $B_j\in \LL(F,E)$ for $j=1,2,...,n$.  
The symbols $\Aa\in \LL(\R^n, \LL(F,G))$ and $\Bb\in \LL(\R^n, \LL(F,E))$ associated to the operators $\A$ and $\B$, are given by 
\begin{align*}
\Aa=\sum_{j=1}^nA_j\otimes e_j,\,\,\, \Bb=\sum_{j=1}^nB_j\otimes e_j\quad \text{and} \quad \Aa(\xi)=\sum_{j=1}^nA_j\xi_j, \,\,\, \Bb(\xi)=\sum_{j=1}^nB_j\xi_j,\quad \xi \in \R^n.
\end{align*}
Define for each $1\leq i,j\leq n$ the matrix
\begin{align}\label{eq:M1}
M_{ij}:=A_i^\ast A_j+B_j^\ast B_i\in \LL(F).
\end{align}
and {\em Laplace coefficients} $\mathbb{M}:\LL(\R^n,F)\to \LL(\R^n,F)$ through 
\begin{align}\label{eq:M2}
\mathbb{M}(X)=\sum_{i,j=1}^n((M_{ij}X(e_j))\otimes e_i.
\end{align}
By Lemma~\ref{lem:Msym}, $\mathbb{M}$ is always self-adjoint. 
If $\mathbb{M}$ is positive semi-definite, we will associate to $\mathbb{M}$ a   first order homogeneous constant coefficient differential operator according to 
\begin{align}\label{eq:DB}
\De_Mu(x)=\sqrt{\mathbb{M}}(Du(x)).
\end{align}
Here $Du(x)\in \LL(\R^n,F)$ denotes the Jacobian matrix of $u$.

For a $C^2$-domain $\Om \subset \R^n$, we consider the boundary condition
\begin{align}
\mathbb{B}(\nu(x))u(x)&=0, \quad x\in \dv \Om, \label{eq:bv2}
\end{align}
 for $u\in C^1(\overline{\Om},F)$, where $\nu$ is the outward pointing unit normal. (By symmetry, the roles of $\Aa$ and $\Bb$ can be interchanged.)
Let $N$ be any $C^1$-extension of $\nu$ to a neighbourhood of $\dv \Om$. Associated to the boundary condition \eqref{eq:bv2} is the \emph{generalized Levi matrix} 
\begin{align}\label{eq:Levi}
L_{B}(x)&=\sum_{j=1}^n\mathbb{B}(e_j)^\ast \mathbb{B}(\dv_jN(x)).
\end{align}
With these notions in place, we prove in Theorem \ref{thm:Weitzen} that the Weitzenb\"ock identity 
\begin{align}\label{eq:Weitzenintro}
\int_{\Om}(\vert \A u(x)\vert^2+\vert \B u(x)\vert^2 -\langle Du(x),\mathbb{M}Du(x)\rangle) dx=\int_{\dv \Om}\langle u(x),L_B(x)u(x)\rangle d\sigma(x)
\end{align}
holds for all $u\in C^1(\overline{\Om},F)$ satisfying \eqref{eq:bv2} on any bounded $C^2$-domain $\Om$. Here $d\sigma$ denotes the surface measure of $\dv \Om$, and $$\langle Du(x),\mathbb{M}Du(x)\rangle= \vert \De_M u(x)\vert^2$$ if $\mathbb{M}$ is positive semi-definite.
If $\{\mathbf{e}_1(x),...,\mathbf{e}_{n-1}(x)\}$ is a local ON-eigenframe of the shape operator $\mathcal{S}_{\dv \Om}^x$ and if $\kappa_1(x),....,\kappa_{n-1}(x)$ are the principal curvatures of $\dv \Om$, we show in Proposition \ref{prop:LeviShape} that the Levi-form, that is, the right hand integrand in \eqref{eq:Weitzenintro}, has the eigenframe expression
\begin{align}\label{eq:Levi2}
\langle u(x),L_{B}(x)u(x)\rangle &=\sum_{j=1}^{n-1}\kappa_j(x)\vert \mathbb{B}(\mathbf{e}_j(x))^\ast u(x)\vert^2,
\end{align}
for $u\in F$.

We have limited us to constant coefficient operators in this paper, but it is 
expected that the extension of \eqref{eq:Weitzenintro} to variable coefficient 
operators and manifolds, contains also a curvature term on the left hand side,
as in the classical case
of Dirac operators. Related to this is a null-lagrangian featuring in our proof
of \eqref{eq:Weitzenintro}, which has the same symmetries as the 
Riemann curvature tensor.

From the Weitzenb\"ock identity \eqref{eq:Weitzenintro}, we derive estimates for the pair $(\A,\B)$ under suitable assumptions on $(\A,\B)$. Recall that an operator $\A$ is of \emph{constant rank} if $\text{rank}\,(\Aa(\xi))=\text{dim}(\text{im}(\Aa(\xi)))$ is the same for all $\xi\in \R^n\setminus\{0\}$. To $\A$ and $\B$ we associate the operator $\mathscr{D}$ by defining 
\begin{align}\label{eq:EllipticOp}
\mathscr{D}u(x):=
\begin{bmatrix}
\A u(x) \\ \B u(x)
\end{bmatrix}.
\end{align}
We characterize in Proposition~\ref{prop:Ddecomp} the operators $\De$
appearing from a complex $(\A,\B)$ in this way.
A homogeneous first order constant coefficient operator $\De$ is called \emph{elliptic} if the symbol $\mathbb{D}(\xi)$ is injective for all $\xi\in \R^n\setminus\{0\}$. Similarly, $\De$ is called \emph{$\C$-elliptic} if the symbol $\mathbb{D}(\zeta)$ is injective for all $\zeta\in \C^n\setminus\{0\}$.  The operator $\De$ from \eqref{eq:EllipticOp} is elliptic if and only if 
\begin{align}\label{eq:EC}
\Ker (\Aa(\xi))=\im(\Bb(\xi)^\ast),
\end{align}
for all $\xi\in \R^n\setminus\{0\}$.
Furthermore, it follows from for example \cite[Lem. 4.1, p. 261]{SW} that whenever $\De$ is elliptic, then both $\A$ and $\B$ are constant rank operators. This leads us to our first main result. 

\begin{Thm}\label{thm:1}
Let $(\A,\B)$ be a pair of first order homogeneous constant coefficient partial differential operators. Assume that $\mathbb{M}$ given by \eqref{eq:M2} is positive semi-definite,
and that $\De_M$ given by \eqref{eq:DB} is $\C$-elliptic. Then for any bounded $C^2$-domain $\Om \subset \R^n$, the coercive estimate 
\begin{align}\label{eq:Coercivity1}
\Vert u\Vert_{W^{1,2}}^2\lesssim \int_{\Om}\vert \A u(x)\vert^2+\vert \B u(x)\vert^2 dx+\int_{\Om}\vert u(x)\vert^2 dx
\end{align}
holds for all $u\in C^1(\overline{\Om}, F)$ satisfying the boundary condition \eqref{eq:bv2}. 
\end{Thm}

We prove Theorem~\ref{thm:1} and compare to coercive estimates in the literature, at the end of Section~\ref{sec:coercive}.
A basic example when Theorem \ref{thm:1} holds is when $E=\Lambda^{k-1}\R^n$, $F=\Lambda^{k}\R^n$, $G=\Lambda^{k+1}\R^n$, $\A=d$, the exterior derivative on differential forms, and $\B=d^\ast=\delta$, the interior derivative on differential forms. 
However, many elliptic operators $\De$ with the boundary conditions \eqref{eq:bv2} imposed are not coercive over $W^{1,2}$, that is, satisfies the estimate \eqref{eq:Coercivity1}. 
J. Kohn and L. Nirenberg showed in \cite{KN} that for pairs of first order differential operators $(\A,\B)$ where $\De$ is elliptic, the solvability of the equation 
\begin{align*}
\B^\ast u(x)=v(x)
\end{align*}
on a domain $\Om\subset \R^n$, where $v(x)$ solves the compatibility condition 
\begin{align*}
\A v(x)=0,
\end{align*}
could still be determined, provided a Morrey inequality of the form 
\begin{align}\label{eq:Morrey}
\int_{\dv \Om} \vert u(x)\vert^2d\sigma(x)\lesssim \int_{\Om}\vert \A u(x)\vert^2+\vert \B u(x)\vert^2 dx+\int_{\Om}\vert u(x)\vert^2dx 
\end{align}
holds for all $u\in C^2(\overline{\Om},F)$ such that $\Bb(\nu)u=0$ on $\dv \Om$. However, no sufficient conditions for \eqref{eq:Morrey} to hold were given in \cite{KN}. 
The inequality \eqref{eq:Morrey} is modelled on the classical Morrey inequality for the Dolbeault complex.
See for example \cite[Sec. 12]{Taylor}.
Extending classical terminology from several complex variables, we say that $\Om$ is \emph{$\B$-strongly pseudoconvex} if the associated Levi form is positive definite
under the boundary condition~\ref{eq:bv2}. 
See Definition~\ref{def:Bsconvex}.
This leads us to our second main result.

\begin{Thm}\label{thm:2}
Let $(\A,\B)$ be a pair of first order homogeneous constant coefficient partial differential operators. Assume that $\mathbb{M}$ is positive semi-definite. Let $\Om\subset \R^n$ be a bounded $C^2$-domain that is $\B$-strongly pseudoconvex. Then the generalized Morrey inequality
\begin{align}\label{eq:Morrey2}
\int_{\dv \Om}\vert u(x)\vert^2d\sigma(x)\lesssim \int_{\Om}(\vert \A u(x)\vert^2+\vert \B u(x)\vert^2) dx
\end{align}
holds
for all $u\in C^1(\overline{\Om},F)$ satisfying the boundary condition \eqref{eq:bv2}. 

As a consequence of the generalized Morrey inequality \eqref{eq:Morrey2},
assuming also that the associated operator $\De$ is elliptic,
the subelliptic square function estimate
 \begin{align}   \label{eq:subE}
  \left(\int_\Om |Du(x)|^2\lambda(x,\dv\Om) dx\right)^{1/2} \lesssim\int_{\Om}(\vert \A u(x)\vert^2+\vert \B u(x)\vert^2) dx
 \end{align}
 holds,
where $\lambda(x,\dv\Om)$ denotes the distance to $\dv\Om$, holds for all $u\in C^1(\overline{\Om},F)$ satisfying the boundary condition \eqref{eq:bv2}.
\end{Thm}

We also show that the left hand side square function norm in \eqref{eq:subE}
dominates $\Vert u\Vert_{W^{1/2,2}(\Om,F)}$, so in particular 
\eqref{eq:subE} implies the usual $1/2$-Sobolev subelliptic estimate. 

Note that the absence of a zero order term $\|u\|_{L^2(\Om)}$ on the right hand
side in \eqref{eq:Morrey2} shows that the hypothesis of Theorem~\ref{thm:2}
implies in particular that the coholomogy vanishes, that is, there are no 
non-trivial $u$ solving $\A u=\B u=0$ in $\Om$ and $\B(\nu)u=0$ on $\dv\Om$.

We further show in Example~\ref{ex:Dolbeaut}, how Theorem~\ref{thm:2} contains
the classical Morrey estimate for the Dolbeaut complex as a special case.

It is natural to ask for a condition on the operator $\B$ such that the class of $\B$-strongly pseudoconvex domains is non-empty and in particular contains all $C^2$ strictly convex domains. Our Proposition \ref{prop:CoCancel} shows that this is true if and only if the operator $\B$ is \emph{cocanceling}. See Definition~\ref{def:CoCancel}. The concept of cocanceling operators was introduced by Van Schaftingen in \cite{VanS} in the study of limiting Sobolev inequalities. It was also an essential condition implicit in the work \cite{BB} of  
Bourgain and Brezis on Hodge systems. It appears that our generalized Levi form together with the cocanceling condition is the first general explicit geometric condition in the literature yielding a Morrey type inequality.

To summarize this paper, we prove sufficient algebraic conditions for coerciveness
for the first order operators $(\A,\B)$ under the ``$1/2$-boundary condition''
\eqref{eq:bv2}, based on the Laplace coefficients $\mathbb{M}$.
\begin{itemize}
\item
$W^{1,2}_0(\Omega)$-coerciveness of $(\A,\B)$, that is, 
\eqref{eq:Coercivity1} for all $u\in W^{1,2}_0(\Omega)$,
 holds whenever $\De$ is elliptic,
as a consequence of Fourier theory.
\item
Full $W^{1,2}(\Omega)$-coerciveness of $(\A,\B)$, that is, 
\eqref{eq:Coercivity1} for all $u\in W^{1,2}(\Omega)$,
 holds whenever $\De$ is $\C$-elliptic,
as a consequence of the Aronszajn-Ne{\u c}as-Smith inequality
Theorem~\ref{thm:ANS}.
\item
$W^{1,2}_B(\Omega)$-coerciveness of $(\A,\B)$, that is, 
\eqref{eq:Coercivity1} for all $u\in W^{1,2}(\Omega)$ satisfying \eqref{eq:bv2},
 holds whenever $\mathbb{M}$ is positive semi-definite and $\De_M$ is $\C$-elliptic.
\item
In particular, $\De_M$ is $\C$-elliptic if $\mathbb{M}$ is positive definite,
as is immediate from \eqref{eq:DB}.
\item
If $\mathbb{M}$ is positive semi-definite, then $|\mathbb{D}_M(\xi)u|= |\mathbb{D}(\xi)u|$ for all $\xi\in \R^n$. See \eqref{eq:Sweeneyident}. 
In particular $\De_M$ is elliptic if and only if $\De$ is elliptic.
Such equivalence is not true for $\C$-ellipticity.
\item
Finally, 
$W^{1/2,2}_B(\Omega)$-coerciveness of $(\A,\B)$, that is, 
\eqref{eq:subE} for all $u\in W^{1,2}(\Omega)$ satisfying \eqref{eq:bv2},
 holds whenever $\mathbb{M}$ is positive semi-definite and $\Om$ is strongly 
 $\B$-pseudoconvex.
\item
And of course, $W^{1,2}(\Omega)$-coerciveness implies $W^{1,2}_B(\Omega)$-coerciveness, which in turn implies 
$W^{1,2}_0(\Omega)$-coerciveness as well as $W^{1/2,2}_B(\Omega)$-coerciveness.
\end{itemize}

It is unclear to us to what extent the hypothesis that $\mathbb{M}$ be 
positive semi-definite in Theorems~\ref{thm:1} and \ref{thm:2} can be relaxed.
In some examples, see end of Section~\ref{sec:coercive}, 
this part of the hypothesis fails, but neverthelss the conclusion holds.
To some extent this is remedied by the fact that one may premultiply 
$\A$ and $\B$ by some invertible matrices, or equivalently change to 
some equivalent norms on $E$, $F$ and $G$.

%=================NEW SECTION=====================================

\subsection{\sffamily  Notation}\label{sec:notation}

Throughout this
paper, we use the notation $X\approx Y$ and $X\lesssim Y$ for estimates to mean that there
exists a constant $C>0$, independent of the variables in the estimate, such that
$X/C\leq Y \leq CX$ and $X\leq CY$, respectively. When $X,Y\in \LL(E,F)$ are linear maps we denote
\begin{align*}
\langle X,Y\rangle:=\mathfrak{R}e[\text{tr}(X^\ast Y)]. 
\end{align*}

%=================NEW SECTION=====================================

\section{\sffamily  Weitzenb\"ock identity}\label{sec:Weitzen}

The goal in this section, is to prove the identity
\eqref{eq:Weitzenintro}, which is fundamental to this paper.

\begin{Def}
The generalized Hodge-Laplacian associated to the pair $(\A,\B)$ is
\begin{align}
\Delta_{\mathscr{D}}&:=\De^\ast  \De= \mathscr{A}^\ast \mathscr{A}+ \mathscr{B}^\ast \mathscr{B}.
\end{align}
\end{Def}

\begin{Lem}\label{lem:HL1}
For $u\in C^2(\Om,F)$, we have
\begin{align}
\Delta_{\mathscr{D}}u(x)&=\Div\mathbb{M}Du(x).
\end{align}
\end{Lem}

\begin{proof}
\begin{align*}
\Delta_{\mathscr{D}}u(x)&=\sum_{1\leq i,j\leq n} (A_i^\ast A_j+B_i^\ast B_j)\dv_i\dv_ju(x)=\sum_{1\leq i,j\leq n} (A_i^\ast A_j+B_j^\ast B_i)\dv_i\dv_ju(x)\\
&=\sum_{1\leq i,j\leq n} M_{ij}\dv_i\dv_ju(x)=\sum_{i=1}^n\dv_i\bigg(\sum_{j=1}^nM_{ij}\dv_ju(x)\bigg)\\
&=\sum_{i=1}^n\dv_i(\sum_{j=1}^n(M_{ij}\dv_ju(x))\otimes e_i)(e_i)=\Div\mathbb{M}Du(x).
\end{align*}
\end{proof}

\begin{Thm}[Weitzenb\"ock identity]
\label{thm:Weitzen}
Let $\Om\subset \R^n$ be a bounded $C^2$ domain and let $\A$ and $\B$ be first order homogeneous constant coefficient operators. 
Then the Weitzenb\"ock identity 
\begin{align}\label{eq:Weitzen}
\int_{\Om}(\vert \A u(x)\vert^2+\vert \B u(x)\vert^2)dx-\int_{\Om}\langle Du(x),\mathbb{M}Du(x)\rangle dx=\int_{\dv \Om}\langle u(x), L_B(x)u(x)\rangle d\sigma(x)
\end{align}
holds for all $u\in C^1(\overline{\Om},F)$ satisfying the boundary condition \eqref{eq:bv2},
where $\mathbb{M}$ is as \eqref{eq:M2} and $L_B(x)$ is as in \eqref{eq:Levi}. 
\end{Thm}

\begin{rem}
We emphasize that identity \eqref{eq:Weitzen} holds for \emph{any} pair of first order homogeneous constant coefficient operators $\A$ and $\B$, and we do not need to assume that the associated operator $\De$ is elliptic.  
\end{rem}

\begin{proof}
We can assume that $u\in C^2(\overline{\Om},F)$, since the general case follows by a limiting argument as in the proof of 
\cite[Thm. 10.3.6]{Ros}.
The basic idea is to apply Stokes' theorem, see for example \cite[Thm. 7.3.9]{Ros}, to the $v$-linear form
\begin{align*}
\theta(x,v)=\langle \Aa(v)u(x),\A u(x)\rangle +\langle\Bb(v)u(x),\B u(x) \rangle -\langle u(x)\otimes v,\mathbb{M}Du(x)\rangle.
\end{align*}
The classical argument is to apply the divergence theorem to the
vector field 
$V(x)=\sum_{j=1}^n\theta(x,e_j)e_j$.
Calculating, using the nabla symbol $\nabla=\sum_{j=1}^ne_j\dv_j$ and
Lemma~\ref{lem:HL1}, we have
\begin{align*}
\text{div}\, V(x)&=\langle \Aa(\nabla )u(\dot{x}),\A u(x)\rangle +\langle \Aa(\nabla )u(x),\A u(\dot{x})\rangle +\langle\Bb(\nabla )u(\dot{x}),\B u(x) \rangle+\langle\Bb(\nabla )u(x),\B u(\dot{x}) \rangle\\
& -\langle u(\dot{x})\otimes \nabla,\mathbb{M}Du(x)\rangle -\langle u(x)\otimes \nabla,\mathbb{M}Du(\dot{x})\rangle\\
&=\vert \A u(x)\vert^2+\langle u(x),\A^\ast \A u(x)\rangle +\vert \B u(x)\vert^2+\langle u(x),\B^\ast \B u(x) \rangle\\
& -\langle Du(x),\mathbb{M}Du(x)\rangle -\langle u(x),\Div \mathbb{M}Du(x)\rangle\\
&=\vert \A u(x)\vert^2+\vert \B u(x)\vert^2 -\langle Du(x),\mathbb{M}Du(x)\rangle.
\end{align*}
Here $\dot{x}$ indicates where the partial derivatives of $\nabla$ act. The divergence theorem and the boundary conditions \eqref{eq:bv2} implies 
\begin{align*}
&\int_{\Om}(\vert \mathscr{A}u(x)\vert^2+\vert \mathscr{B}u(x)\vert^2-\langle Du(x),\mathbb{M}Du(x)\rangle)dx
=\int_{\dv\Om}\theta(x,\nu(x)) d\sigma(x)\\
&=\int_{\dv\Om}\langle u(x),\mathbb{A}(\nu(x))^\ast \mathscr{A}u(x)-(\mathbb{M}Du(x))(\nu(x))\rangle d\sigma(x).
\end{align*}
Furthermore
\begin{align*}
(\mathbb{M}Du(x))(\nu(x))&=\sum_{i,j}^nM_{ij}\nu_i(x)\dv_ju(x)=\sum_{i,j}^n(A_i^\ast A_j+B_j^\ast B_i)\nu_i(x)\dv_ju(x),
\end{align*}
and so 
\begin{align*}
&\mathbb{A}(\nu(x))^\ast\mathscr{A}u(x)-(\mathbb{M}Du(x))(\nu(x))\\
&=\sum_{i,j=1}^n(A_i^\ast A_j-A_i^\ast A_j-B_j^\ast B_i)\nu_i(x)\dv_ju(x)\\
&=-\sum_{i,j=1}^nB_j^\ast B_i\nu_i(x)\dv_ju(x)=-\sum_{j=1}^nB_j^\ast \mathbb{B}(\nu(x))\dv_ju(x).
\end{align*}
Note that it is at this stage in the computation that the form of $\mathbb{M}$ is determined by the condition that we want the $A_i^\ast A_j$
 terms to cancel since we have no knowledge of $\Aa(\nu(x))u(x)$ on $\dv \Om$.  Using a $C^1$ extension $N$ of $\nu$ in a neighbourhood of $\dv \Om$, we calculate

\begin{align*}
\sum_{j=1}^nB_j^\ast \mathbb{B}(\nu(x))\dv_ju(x)&=\sum_{j=1}^nB_j^\ast \dv_j\big(\mathbb{B}(N(x))u(x)\big)-\sum_{j=1}^n\big(B_j^\ast \dv_j\mathbb{B}( N(x))\big)u(x)\\
&=-\mathscr{B}^\ast\big(\mathbb{B}(N(x))u(x)\big)+\big(\mathscr{B}^\ast\mathbb{B}( N(x))\big)u(x).
\end{align*}

We now write $\mathscr{B}^\ast$ in a way that is adapted to the boundary $\dv \Om$. Let $\{\mathbf{e}_j(x)\}_{j=1}^{n-1}$ be a local ON-frame for $T\dv \Om$ around a point $x\in \dv \Om$. The nabla symbol can be written 
\begin{align} \label{eq:nablasplit}
\nabla=\sum_{j=1}^ne_j\dv_j=\nu(x)\dv_{\nu(x)}+\sum_{j=1}^{n-1}\mathbf{e}_j(x)\dv_{\mathbf{e}_j(x)}=\nu(x)\dv_{\nu(x)}+\nabla^T. 
\end{align}
Thus
$\mathscr{B}^\ast =-\mathbb{B}(\nabla)^\ast=-\mathbb{B}(\nu(x))^\ast \dv_\nu -\mathbb{B}(\nabla^T)^\ast$.

This gives
\begin{align*}
\langle u(x),-\mathscr{B}^\ast\big(\mathbb{B}(N(x))u(x)\big)\rangle&=\langle u(x),\mathbb{B}(\nu(x))^\ast\dv_\nu\big(\mathbb{B}(N(x))u(x)\big)\rangle+\langle u(x),\mathbb{B}(\nabla^T)^\ast\big(\mathbb{B}(\N(x))u(x)\big)\rangle\\
&=\langle \mathbb{B}(\nu(x))u(x),\dv_\nu\big(\mathbb{B}(N(x))u(x)\big)\rangle+\langle u(x),\mathbb{B}(\nabla^T)^\ast\big(\mathbb{B}(\nu(x))u(x)\big)\rangle
=0,
\end{align*}
since $ \mathbb{B}(\nu(x))u(x)=0$ and $\mathbb{B}(\nabla^T)^\ast\big(\mathbb{B}(\nu(x))u(x)\big)=0$ for $x\in \dv \Om$. Finally, 
\begin{align*}
-\mathscr{B}^\ast\mathbb{B}( N(x))=\sum_{j=1}^n\mathbb{B}(e_j)^\ast\mathbb{B}(\dv_jN(x))=L_B(x)
\end{align*}

Thus, $\langle u(x),\mathbb{A}(\nu(x))^\ast \mathscr{A}u(x)-(\mathbb{M}Du(x))(\nu(x))\rangle=\langle u(x),L_B(x)u(x)\rangle$, which concludes the proof. 

\end{proof}

\begin{rem}\label{remNull}
The proof of Theorem~\ref{thm:Weitzen} shows in particular that 
$\vert \A u(x)\vert^2+\vert \B u(x)\vert^2 -\langle Du(x),\mathbb{M}Du(x)\rangle$ is a quadratic \emph{null-lagrangian}. 
Quadratic null-lagrangians are precisely those quadratic forms $Q_N(X)=\langle X, \mathbf{N}X\rangle$, where $X\in \LL(\R^n,F)$ and $\mathbf{N}\in \LL(\LL(\R^n,F),\LL(\R^n,F))$ such that $Q_N(u\otimes v)=0$ for all $u\in F$ and $v\in \R^n$. 
More generally, for definitions and properties of null-lagrangians we refer the reader to \cite{BCO,OS}. 
\end{rem}

We include in this section a number of properties of auxiliary results that will be used in later sections, in particular about the Laplace coefficients 
$\mathbb{M}$. 

\begin{Lem}  \label{lem:Msym}
Let $\{v_1,...,v_N\}$ be an ON-basis for $F$ and $\{e_1,...,e_n\}$ be an ON-basis for $\R^n$. Then 
\begin{align*}
\mathbb{M}^{\alpha \beta}_{ij}=\langle v_\alpha \otimes e_i,\mathbb{M}(v_\beta\otimes e_j)\rangle=\langle A_iv_\alpha, A_j v_\beta \rangle+\langle B_jv_\alpha, B_i v_\beta\rangle.
\end{align*}
In particular $\mathbb{M}^\ast=\mathbb{M}$.
\end{Lem}

\begin{proof}
The Laplace coefficients are computed as
\begin{align*}
&\mathbb{M}^{\alpha \beta}_{ij}=\langle v_\alpha \otimes e_i,\mathbb{M}(v_\beta\otimes e_j)\rangle=\sum_{k,l=1}^n\langle v_\alpha \otimes e_i,((M_{kl}(v_\beta\otimes e_j)(e_k))\otimes e_l\rangle\\
&=\sum_{k,l=1}^n\langle (v_\alpha \otimes e_i)(e_l),M_{kl}v_\beta \delta_{jk}\rangle=\sum_{k,l=1}^n\langle v_\alpha,M_{kl}v_\beta \rangle \delta_{jk}\delta_{il}\\
&=\langle v_\alpha,M_{ji}v_\beta \rangle=\langle M_{ij}v_\alpha,v_\beta \rangle=\langle A_iv_\alpha, A_j v_\beta \rangle+\langle B_jv_\alpha, B_i v_\beta\rangle.
\end{align*}
In particular $\mathbb{M}^{\beta\alpha}_{ji}= \mathbb{M}^{\alpha \beta}_{ij}$, which completes the proof.
\end{proof}

The Laplace coefficients $\mathbb{M}$ need not be positive semi-definite, even 
if the associated operator $\De$ is elliptic, as the following example shows.
 
\begin{ex}\label{ex:1}
  For $\eps\in\R$, consider
\begin{align*}
\xymatrix{C^\infty(\Om,\Lambda^0\R^2) & \ar[l]_{\delta} C^\infty(\Om,\Lambda^1\R^2)\ar[r]^{\eps d} &C^\infty(\Om,\Lambda^2\R^2)},
\end{align*}
where $\B=\delta$ is divergence and $\A=\eps d$ is curl, rescaled by $\eps$.
When $\epsilon=1$, the anticommutation relation, 
see \cite[Thm. 2.8.1]{Ros} and the notation for interior and exterior products used there, shows that
$$
  M_{ij}u=e_i\ri (e_j\wedg u)+ e_j\wedg(e_i\ri u)= \delta_{ij}u, \qquad i,j=1,2, 
$$ 
and $\mathbb{M}$ is the identity map.
However when $\eps>0$ is small, the map $\mathbb{M}$ is indefinite.
Indeed, for
$$
X=\begin{bmatrix}
0 & -1\\ 1 & 0
\end{bmatrix}
$$
a direct calculation shows that 
$\langle X,\mathbb{M}X \rangle=-2+4\eps^2$.
\end{ex}

For a given pair $(\A,\B)$ of first order homogeneous constant coefficient partial differential operators, we may replace $\A$ and $\B$ by
$A_0\A$ and $B_0\B$, where $A_0\in \text{GL}(G)$ and $B_0\in \text{GL}(E)$ are two invertible maps.
This does not affect ellipticity of $\De$, nor the boundary condition 
\eqref{eq:bv2}, and clearly
$$
\int_{\Om}(\vert A_0\A u(x)\vert^2+\vert B_0\B u(x)\vert^2)dx
\approx
\int_{\Om}(\vert \A u(x)\vert^2+\vert \B u(x)\vert^2)dx.
$$ 
However, the Laplace coefficients $\mathbb{M}$ are now replaced
$$
  M_{ij}\mapsto A_i^*A_0^*A_0 A_j+ B_j^*B_0^*B_0B_i,
$$
and this may affect positivity of $\mathbb{M}$.
Indeed, the choice $A_0=\epsilon$ and $B_0=1$ in Example~\ref{ex:1}
turn a positive definite $\mathbb{M}$ into an indefinite $\mathbb{M}$.
Therefore, the conclusion of
Theorem \ref{thm:1} and Theorem \ref{thm:2} can be obtained by 
verifying the hypothesis, not for $(\A,\B)$, but for 
some auxiliary $(A_0\A, B_0\B)$.

Even if $\mathbb{M}$ may fail to be positive semi-definite for general
$X$, it is always positive semi-definite on rank one matrices $X$.
We say that $X\mapsto \langle X,\mathbb{M}X \rangle$ 
is {\em rank-one convex} if
\begin{align}\label{eq:rankone}
\langle u\otimes v, \mathbb{M}(u\otimes v)\rangle\geq 0
\end{align}
for all $u\in F$ and all $v\in \R^n$. By \cite[Lem. 5.27, p. 192]{Da}, \eqref{eq:rankone} is equivalent to the standard definition of rank-one convexity given in \cite[Def. 5.1, p. 156]{Da}.

\begin{Lem}\label{lem:LH}
For any $u\in F$ and $v\in \R^n$, we have
\begin{align*}
\langle u\otimes v, \mathbb{M}(u\otimes v)\rangle=\langle u, \mathbb{D}(v)^\ast \mathbb{D}(v) u\rangle= \vert \mathbb{D}(v) u\vert^2,
\end{align*}
where $\mathbb{D}$ is the symbol of $\De$. 
In particular
$X\mapsto \langle X,\mathbb{M}X \rangle$ 
is rank-one convex, for any $(\A,\B)$.

Moreover, the quadratic form $X\mapsto \langle X,\mathbb{M}X \rangle$ satisfies the Legendre-Hadamard condition 
\begin{align*}
\langle u\otimes v, \mathbb{M}(u\otimes v)\rangle\geq \lambda  \vert v\vert^2\vert u\vert^2
\end{align*}
for some $\lambda>0$ if and only if the associated operator $\De$ is elliptic. 
\end{Lem}

\begin{proof}
Calculating, we have
\begin{align*}
\langle u\otimes v, \mathbb{M}(u\otimes v)\rangle&=\sum_{i,j=1}^n\langle u\otimes v, (M_{ij}(u\otimes v)(e_j))\otimes e_i\rangle=\sum_{i,j=1}^n\langle u\langle v,e_i\rangle, M_{ij}u\langle v,e_j\rangle \rangle\\
&=\sum_{i,j=1}^nv_iv_j\langle u, M_{ij}u\rangle=\sum_{i,j=1}^nv_iv_j\langle u ,(A_i^\ast A_j+B_j^\ast B_i)u\rangle\\
&=\sum_{i,j=1}^nv_iv_j(\langle A_iu ,A_ju\rangle+\langle B_ju ,B_iu\rangle)=\sum_{i,j=1}^n(\langle A_iv_iu ,A_jv_ju\rangle+\langle B_jv_ju ,B_iv_iu\rangle)\\
&=\langle \Aa(v)u, \Aa(v)u\rangle+\langle\Bb(v)u,\Bb(v)u \rangle=\langle \mathbb{D}(v)u,\mathbb{D}(v)u\rangle=\langle u, \mathbb{D}(v)^\ast \mathbb{D}(v) u\rangle,
\end{align*}
from which rank one convexity, as well as the equivalence of the Legendre-Hadamard condition and ellipticity, is immediate.
\end{proof}

We end this section by considering the converse of \eqref{eq:EllipticOp}.
Namely, suppose that we are given a first order homogeneous constant coefficient operator $\De:C^\infty(\R^n,F)\to C^\infty(\R^n,H)$. 
When does $\De$ admit a decomposition $\De=\A\oplus \B$, so that we
can apply the Weitzenb\"ock identity
in Theorem~\ref{thm:Weitzen} to $(\A, \B)$?
Note that a decomposition $\De=\A\oplus \B$ of a given first order elliptic operator $\De$ brings with it the local boundary condition \eqref{eq:bv2}.
Without such a decomposition of $\De$, it is not clear what local boundary
conditions that are natural for $\De$.
In most applications, $(\A, \B)$ are associated to a short complex
\begin{align}   \label{eq:complex}
\xymatrix{C^\infty(\Om,E)  \ar[r]^{\B^*} & C^\infty(\Om,F)\ar[r]^{\A} &C^\infty(\Om,G)},
\end{align}
meaning that the $\A\circ \B^\ast=0$.

\begin{Prop} \label{prop:Ddecomp}
Let $\De: C^\infty(\R^n,F)\to C^\infty(\R^n,H)$ be a 
first order homogeneous constant coefficient operator. 
Consider orthogonal splittings $H= G\oplus E$, with $\mathbf{P}$ denoting 
orthogonal projection onto $E$. 
Define 
$\A=(I-\mathbf{P})\De: C^\infty(\R^n,F)\to C^\infty(\R^n,G)$
and $\B=\mathbf{P}\De : C^\infty(\R^n,F)\to C^\infty(\R^n,E)$.
Then $\De=\A\oplus \B$, under suitable identifications.

Moreover, $\A\circ \B^\ast=0$ holds if and if $E$ is an invariant subspace
of $\Di(\xi)\circ \Di(\xi)^\ast$ for all $\xi \in S^{n-1}$.
\end{Prop}

\begin{proof}
Only the last statement needs to be proved. 
But
\begin{align*}
\Di(\xi)\circ \Di(\xi)^\ast&=\begin{bmatrix}
\Aa(\xi) \\ \Bb(\xi)
\end{bmatrix}
\begin{bmatrix}
\Aa(\xi)^\ast & \Bb(\xi)^\ast
\end{bmatrix}
=\begin{bmatrix}
\Aa(\xi)\Aa(\xi)^\ast &\Aa(\xi)\Bb(\xi)^\ast\\
\Bb(\xi) \Aa(\xi)^\ast & \Bb(\xi)\Bb(\xi)^\ast
\end{bmatrix},
\end{align*}
in the splitting $H= G\oplus E$.
Since $E$ is an invariant subspace if and and only if the (1,2) block
in the matrix vanishes, the proof is complete.
\end{proof}

\begin{ex}
Let $\De u(x)=\text{Def}\,u(x)=\frac{1}{2}(Du(x)+Du(x)^\ast)$ be the symmetric gradient, which is an overdetermined elliptic operator.  
In $\R^2$, the symmetric gradient is equivalent to the operator with symbol
$
\mathbb{D}(\xi)=\begin{bmatrix}
\xi_1 & 0 \\ \tfrac 12\xi_2 & \tfrac 12\xi_1 \\ 0 & \xi_2
\end{bmatrix}.
$
Then
\begin{align*}
\mathbb{D}(\xi)\circ \mathbb{D}(\xi)^\ast=
\frac{1}{4}\begin{bmatrix}
4\xi_1^2 & 2\xi_1\xi_2  & 0\\
2\xi_1\xi_2 &\xi_1^2+\xi_2^2&  2\xi_1\xi_2\\
0 &2\xi_1\xi_2 & 4\xi_2^2
\end{bmatrix}
\end{align*}
This has no non-trivial subspaces invariant for all $\xi$.
Indeed, the choices $\xi=(1,0)$ and $\xi=(0,1)$ reveals that any non-zero such
subspace must contain $(1,0,0)$, $(0,1,0)$ or $(0,0,1)$. Varying $\xi$ then implies that
the subspace must contain all vectors. Thus the symmetric gradient has no non-trivial decomposition into operators $\A$ and $\B$ besides the trivial choices $\A=0$ and $\B=\text{Def}$, or vice versa. 
\end{ex}

\begin{ex}
Let $\De: C^\infty(\R^n, E)\to C^\infty(\R^n, F)$ be a \emph{Dirac type} operator. This means in particular that 
\begin{align*}
\mathbb{D}(\xi)\circ \mathbb{D}(\xi)^\ast=\vert \xi\vert^2\text{id}_F,
\end{align*}
which is a diagonal operator. Therefore any subspace $E$ is invariant for all $\mathbb{D}(\xi)$. 
\end{ex}

%========================NEW SECTION==========================================

\section{\sffamily  Coerciveness and the ANS inequality}   \label{sec:coercive}

For the de Rham complex with $\A=d$ and $\B=\delta$, we have the classical Gaffney inequality
\begin{align*}
\int_{\Om}\vert Du(x)\vert^2dx \lesssim\int_{\Om}(\vert du(x)\vert^2+\vert \delta u(x)\vert^2+\vert u(x)\vert^2)dx
\end{align*}
on a $C^2$ domain $\Om$, for all $k$-forms $u\in C^1(\overline{\Om},\Lambda^k \R^n)$ that are tangential on $\dv \Om$, that is, $\nu(x)\ri u(x)=0$, $x\in \dv \Om$.
In this section, we prove Theorem~\ref{thm:1}, which gives the analogous 
Gaffney inequality for a class of first order homogeneous constant coefficient operators $(\A,\B)$.
In the special case when $\mathbb{M}$ is positive definite, the coercive 
estimate \eqref{eq:Coercivity1} is an immediate corollary of Theorem~\ref{thm:Weitzen}.
For example, $\mathbb{M}$ is the identity for the de Rham complex above.
More generally, as discussed in Section~\ref{sec:Weitzen}, we may
also deduce the coercive estimate by checking that $\mathbb{M}$
corresponding to $(A_0\A, B_0 \B)$, for some invertible matrices $A_0, B_0$,
is positive definite.

The proof of Theorem~\ref{thm:1}
in the case when $\mathbb{M}$ is only positive semi-definite is more
subtle, and uses the following coercive estimate by
Aronszajn, Ne{\u c}as and Smith (ANS).

\begin{Thm}[Aronszajn-Ne{\u c}as-Smith]
\label{thm:ANS}
Let $\Om\subset \R^n$ be a bounded $C^2$ domain and let $\Pop: C^\infty(\R^n, E)\to C^\infty(\R^n, F)$ be a first order homogeneous constant coefficient operator.  
Then 
\begin{align}\label{eq:ANS}
\Vert u\Vert_{W^{1,2}(\Om,E)}\lesssim \Vert \Pop u\Vert_{L^2(\Om,F)}+\Vert u\Vert_{L^2(\Om,E)}
\end{align}
for all $u\in W^{1,2}(\Om,E)$,
if and only if
the symbol $\Pp$ of $\Pop$ is $\C$-elliptic, that is, $\Pp(\zeta)$ is injective for all $\zeta\in\C^n\setminus\{0\}$.
\end{Thm}

Note carefully that the ANS-inequality holds \emph{without} imposing \emph{any} boundary conditions on $u$, as a consequency of the $\C$-ellipticity, which is a purely algebraic necessary and sufficient condition.
Verifying the condition of $\C$-ellipticity is in principle much simpler than verifying either the necessary and sufficient conditions in the general case in \cite[Thm. 1]{Sw1} and \cite{DeF}.

We remark that the ANS-inequality holds in more general situations: $\Pop$ may be a higher order differential operator, with variable coefficients, and the principal symbol need only be $\C$-elliptic at the boundary. Moreover, $\Om$
may be a bounded Lipschitz domain (unlike for Theorem~\ref{thm:1}) and $L^2$ Sobolev norms may be replaced by $L^p$ Sobolev norms, $1<p<\infty$.

Since the ANS-inequality is not as well-known as it should be,
we include a proof of Theorem~\ref{thm:ANS} for the convenience of the reader. This follows \cite{Sm61} rather than the later proof in \cite{Sm70}. 

\begin{proof}
We first prove the necessity of $\C$-ellipticity, following \cite{Nec}: assume that there exists
$v\in \C^n$, $v\neq 0$, and $\zeta\in\C^n\setminus\{0\}$, such that $\Pp(\zeta)v=0$. For each $\lambda\in \C$, $\lambda\neq 0$,  consider the plane wave
\begin{align*}
u_{\lambda}(x)=\lambda^{-1}ve^{\lambda \langle \zeta,x\rangle}.
\end{align*}
Then 
\begin{align*}
\Pop u_\lambda(x)=\sum_{j=1}^n\lambda^{-1}P_jv\dv_je^{\lambda \langle \zeta,x\rangle}=\sum_{j=1}^nP_j\zeta_jve^{\lambda \langle \zeta,x\rangle}=\Pp(\zeta)ve^{\lambda \langle \zeta,x\rangle}=0.
\end{align*}
Thus $u_\lambda$, is a solution and so is any linear combination of such plane wave solutions. 
Let $V$ be the $W^{1,2}(\Om)$ closed linear span of $\{u_\lambda\}_{\lambda \in \C}$, which is an infinite dimensional Hilbert space.
If \eqref{eq:ANS} holds, then the $L^2$ and $W^{1,2}$ norms on $V$ are comparable.
This contradicts the Rellich-Kondrachov compactness theorem.

We next prove the sufficiency of $\C$-ellipticity, following \cite{Sm61}.
We aim to prove that there exists $k<\infty$, such that for each multi-index $\alpha=(\alpha_1,\ldots,\alpha_n)$ with $\vert \alpha \vert>n(k-1)$, there exists a symbol $\Qq_\alpha(\xi)\in \LL(F,E)$ homogeneous of degree $|\alpha|-1$, such that
\begin{align}\label{eq:AlgeComp}
\mathbb{Q}_\alpha(\xi) \mathbb{P}(\xi)= \xi^\alpha\text{id}_E. 
\end{align}
This will complete the proof, since then $\dv^\alpha u(x)=\mathscr{Q}_\alpha (\mathscr{P}u(x))$, and either J.L. Lions lemma \cite[Thm. 6.11-5, p. 381]{Cia} or a Calder\'on-Zygmund type argument as in \cite{Sm61}
applies to prove \eqref{eq:ANS}.

We may assume that $E=\C^m$ and $F=\C^N$.
To find $k$ and construct $\Qq_\alpha(\zeta)$, we consider 
$m\times m$-minors $\Pp_s(\zeta)$ of $\Pp(\zeta)$ formed by selecting rows $i_1, i_2, \ldots, i_m$ in $\Pp(\zeta)$, where $s=\{i_1, i_2, \ldots, i_m\}\subset \{1, 2, \ldots, N\}$.
Define the subdeterminant $p_s(\zeta)=\det(\Pp_s(\zeta))$. 
By the Cauchy-Binet formula 
\begin{align*}
\det(\mathbb{P}(\zeta)^\ast \circ \mathbb{P}(\zeta))=\sum_{\vert s\vert=m}\vert p_s(\zeta)\vert^2.
\end{align*}
Thus, $\Pp(\zeta)$ fails to be injective if and only if the scalar polynomials $\{p_s(\zeta)\}_s$ have a common zero, for any $\zeta \in \C^n\setminus \{0\}$. 
Assuming $\C$-ellipticity, we conclude that $\{p_s(\zeta)\}_s$ have a common zero
only for $\zeta=0$. 
The Hilbert-Nullenstellensatz from algebraic geometry 
shows that for each $j$, there exist $k_j\in\Z$ such that
$\zeta_j^{k_j}=\sum_s q_{s,j}(\zeta) p_s(\zeta)$ for some polynomials $q_{s,j}(\zeta)$, since $\zeta_j$ vanishes at $\zeta=0$.

Set $k=\max(k_j)$ and consider any multi-index $\alpha$ with $|\alpha|>n(k-1)$.
By the pigeonhole principle, there exists $j$ such that $\alpha_j\ge k$.
Then $\zeta_j^{k_j}$ divides $\zeta^\alpha$, and it follows that  
$$
\zeta^\alpha=\sum_s q_{s,\alpha}(\zeta) p_s(\zeta)
$$ 
for the polynomials $q_{s,\alpha}(\zeta)= (\zeta^\alpha/\zeta_j^{k_j})  q_{s,j}(\zeta)$.
To construct $\Qq_\alpha(\xi)$, we recall the cofactor formula
\begin{align*}
\text{cof}(\Pp_s(\xi))\Pp_s(\xi)=p_s(\xi)\text{id}_E. 
\end{align*} 
Let $R_s$ be the $m/N$ matrix with entries $0$ and $1$ such that
$\Pp_s(\xi)= R_s \Pp(\xi)$, and set
$$
  \Qq_\alpha(\xi)= \sum_s q_{s,\alpha}(\xi) \text{cof}(\Pp_s(\xi)) R_s.
$$ 
Then it follows that 
$$
  \mathbb{Q}_\alpha(\xi) \mathbb{P}(\xi)=
  \sum_s q_{s,\alpha}(\xi) \text{cof}(\Pp_s(\xi)) R_s \mathbb{P}(\xi)=
  \sum_s q_{s,\alpha}(\xi) p_s(\xi)\text{id}_E =  
   \xi^\alpha\text{id}_E,
$$
as needed.
\end{proof}

\begin{ex}
The classical example where the ANS-inequality applies, is the symmetric gradient
$\text{Def}\,u(x)=(Du(x)+Du(x)^\ast)/2$.
In $\R^2$, the symmetric gradient is equivalent to the operator
\begin{align*}
\Pop u(x)=\begin{bmatrix}
1 & 0 \\ 0 &1\\ 0 & 0
\end{bmatrix}
\dv_xu+
\begin{bmatrix}
0 & 0 \\ 1 &0\\ 0 & 1
\end{bmatrix}
\dv_yu.
\end{align*}
Here is computation shows that
\begin{align*}
\det(\Pp(\zeta)^\ast \Pp(\zeta))=(\vert \zeta_1\vert^2+\vert \zeta_2\vert^2)^2-\vert \xi_1\vert^2\vert \xi_2 \vert^2.
\end{align*}
Thus $\Pop$ is $\C$-elliptic, and the ANS-inequality yields the classical Korn inequality.  
\end{ex}

\begin{ex}
For the Cauchy--Riemann system 
\begin{align*}
\mathcal{P}=\begin{bmatrix}
\dv_x & -\dv_y\\ \dv_y & \dv_x
\end{bmatrix}
\end{align*}
in $\R^2=\C$, one has
$\det(\mathbb{P}(\zeta))= \zeta_1^2+ \zeta_2^2$.
Here $\mathbb{P}(\zeta)$ is not injective on the two complex lines $\zeta_2=\pm i\zeta_1$, and it follows that $\mathcal{P}$ is elliptic, but not $\C$-elliptic. 
\end{ex}

\begin{proof}[Proof of Theorem \ref{thm:1}]
From Theorem~\ref{thm:Weitzen}, and the assumptions that $\mathbb{M}$ is positive
semi-definite and that $\De_M$ is $\C$-elliptic, we obtain the estimate
\begin{align*}
\Vert u\Vert_{W^{1,2}(\Om,F)} & \lesssim 
\int_{\Om}(\vert \De_M u(x)\vert^2+\vert u(x)\vert^2)dx \\
& \lesssim\int_{\Om}(\vert \A u(x)\vert^2+ \vert \B u(x)\vert^2+\vert u(x)\vert^2)dx
+ \int_{\dv\Om} \vert u(x)\vert^2dx,
\end{align*}
using also that the Levi-form $L_B(x)$ is bounded for $C^2$ domains.
Noting as in \cite[Cor. 10.3.8]{Ros} the trace estimate
$$
\int_{\dv\Om} \vert u(x)\vert^2dx\lesssim\int_{\Om}(\vert Du(x)\vert \vert u(x)\vert+\vert u(x)\vert^2) dx,
$$
absorption of the $|Du(x)|$-term yields the coercive estimate stated in
Theorem \ref{thm:1}.
\end{proof}

\begin{rem}\label{rem:Thm1}
Note that the associated operator $\mathscr{D}$ in Theorem \ref{thm:1} must necessarily be elliptic. This is due to the fact \eqref{eq:Coercivity1} implies coercivity over $W^{1,2}_0(\Om,E)$, which holds if and only if $\mathbb{M}$ satisfy the Legendre-Hadamard condition, see \cite[Thm. 3.42]{GM}, and hence by Lemma \ref{lem:LH}, $\mathscr{D}$ must be elliptic. 
\end{rem}

We conclude this section by a comparison to related coercive estimates 
in the literature. 
Sweeney \cite[Thm. 1, p. 384]{Sw1} gives a necessary and sufficient condition for the coercive estimate \eqref{eq:Coercivity1} to hold under the boundary condition \eqref{eq:bv2}. This result covers operators with variable coefficients and lower order terms, but the 
condition is not algebraic.
There are also necessary and sufficient conditions for coerciveness by Agmon~\cite{Agmon} and  De Figueiredo~\cite{DeF}. 
Their conditions are similar in flavour to Sweeney's. 

In \cite[Thm. 4, p. 387]{Sw1}, an algebraic condition which is sufficient for coerciveness, is given. 
The standard assumption in \cite{Sw1} is that $\B^*$ and $\A$ form an elliptic complex \eqref{eq:complex}, that is, that the symbol sequence
\begin{align*}
\xymatrix{E  \ar[r]^{\Bb^*(\xi)} & F \ar[r]^{\Aa(\xi)} & G}
\end{align*}
is exact for all $\xi\in \R^n\setminus\{0\}$.
For constant coefficient operators, 
the condition for coerciveness in \cite[Thm. 4]{Sw1} amounts to that this exactness continues to hold for all $\zeta\in \C^n\setminus\{0\}$.
To compare this to the condition for coerciveness given in 
our Theorem \ref{thm:1}, we note by a calculation similar to the 
proof of Lemma~\ref{lem:LH}, that 
\begin{equation}  \label{eq:Sweeneyident}
   \langle u\otimes \zeta, \mathbb{M}(u\otimes \zeta)\rangle=
   \vert \Aa(\zeta)u \vert^2 + \vert \Bb(\overline{\zeta})u \vert^2
   =   \vert \Aa(\zeta)u \vert^2 + \vert (\Bb^*(\zeta))^*u \vert^2
\end{equation}
for any $u\in F$ and $\zeta\in \C^n$.
Thus, assuming that $\mathbb{M}$ is positive semi-definite,
the $\C$-ellipticity conditions in \cite[Thm. 4]{Sw1} and in
Theorem \ref{thm:1} are equivalent.
It is interesting to note that the condition in \cite[Thm. 4]{Sw1}
is symmetric with respect to $\A$ and $\B^*$, in contrast to 
the condition in Theorem~\ref{thm:1}.
 
In some examples, $\B^*$ and $\A$ form an elliptic complex
and the condition in \cite[Thm. 4]{Sw1} holds, but our condition fails since
$\mathbb{M}$ is indefinite. We have seen this in Example~\ref{ex:1}, and 
another is the following.

\begin{ex}[Donaldson-Sullivan complex]
In their work \cite{DS} on quasiconformal 4-manifolds, Donaldson and Sullivan 
use the short elliptic complex
\begin{align*}
\xymatrix{C^\infty(\Om,\Lambda^0\R^4) \ar[r]^{d}& C^\infty(\Om,\Lambda^1\R^4)\ar[r]^{d^+} &C^\infty(\Om,\Lambda^2_+\R^4)}.
\end{align*}
Here 
$$
  \A u= d^+ u := \tfrac 12(1+*) du, 
$$
where $\tfrac 12(1+*)$ is the projection onto the Hodge star self-dual 
bivectors in $\R^4$,
and
$$
 \B u= d^* u
$$
is divergence of the vector field $u$.
In \cite{DS} these maps are considered on closed manifolds, but if one 
considers the coerciveness on bounded domains in $\R^4$ under the
boundary condition \eqref{eq:bv2}, then this follows from the
condition in \cite[Thm. 4]{Sw1}.
It is interesting to note that coerciveness thus holds even if the
projection $\tfrac 12(1+*)$ is not invertible. (Though it is injective on the range of any map $\Aa(\xi)$.)
We also note that the Laplace coefficients for this example $(\A,\B)$
are indefinite. Indeed, $\langle X, \mathbb{M}X \rangle=-1$ for the linear 
vector field whose Jacobian matrix equals
$$
  X= \begin{bmatrix}
    0 & -1 & 0 & 0 \\
    1 & 0 & 0 & 0 \\
    0 & 0 & 0 & 0 \\
    0 & 0 & 1 & 0
  \end{bmatrix}.
$$
\end{ex}

There are of course examples $(\A,\B)$ where $\B^*$ and $\A$ do not
form an elliptic complex.
In this case \cite[Thm. 4]{Sw1} does not apply, but in some examples the
condition in Theorem~\ref{thm:1} can be verified.

\begin{ex}[Symplectic de Rham complex]
Consider the standard de Rham complex
\begin{align*}
\xymatrix{C^\infty(\Om,\Lambda^{k-1}\R^n) \ar[r]^{d} & C^\infty(\Om,\Lambda^k\R^n)\ar[r]^{d} &C^\infty(\Om,\Lambda^{k+1}\R^n)}
\end{align*}
on $\R^n$. 
Let $\A=d$.
Let $J\in\LL(\R^n,\R^n)$ be a non-trivial isometry and consider the 
duality
$$
  \langle u, v\rangle_J= \langle Ju, v \rangle.
$$
With the respect to this duality (which in general is not an inner product), the adjoint of $d$ is
$\B= -J^*\delta J$, 
where we continue to denote the $\wedg$-homomorphism extending $J$,
by $J$.
It follows from \cite[Prop. 2.7.1]{Ros} that its symbol is 
$$
  \Bb(\xi)u= -(J^*\xi)\lctr u. 
$$
Note that $\Bb(\xi)^*u= -(J^*\xi)\wedg u$, and that $\xi\wedg(J^*\xi)\ne 0$
for some $\xi\ne 0$.
For $(\A,\B)$ we therefore note that $\B^*$ and $\A$ do not 
form an elliptic complex.
As already noted, $\mathbb{M}=I$ and Theorem~\ref{thm:1} applies
when $J$ is the identity.
By continuity Theorem~\ref{thm:1} also applies when $J$ is a small rotation.

An example of a large perturbation $J$ of the identity comes from 
symplectic geometry.
The symplectic de Rham complex was introduced by Brylinski in \cite{Bry} in the study of smooth closed symplectic manifolds and has been studied in  subsequent works \cite{Ma,Yan, TY1,TY2,TY3}.
Here $J$ is the standard symplectic structure on $\R^n$, with $n$ even.
However, in this case $(\A,\B)$ is not even elliptic, as it is readily verified that 
$$
\begin{cases}
  \xi\wedg u=0, \\
  (J^*\xi)\lctr u=0,
\end{cases}
$$
has many non-trivial solutions.  In particular, none of the estimates
\eqref{eq:Coercivity1} or \eqref{eq:subE} can hold.
However, restricting $\Lambda^{k}\R^n$ to an eigenspace of $J$ will restore
ellipticity. 
Replacing $\B$ by $J^\ast \B=\delta J$ and using \cite[Prop. 2.7.1]{Ros}, 
the symbol of $J^\ast \B$ equals 
\begin{align*}
J^\ast \Bb(\xi)u=J^\ast((J \xi)\lctr u )= \xi\lctr (J^\ast u). 
\end{align*}
Let $k$ be even and let $\Lambda^k\R^n=V_+^k\oplus V_-^k$ be the eigenspace decomposition with respect to $J$ and note that $J^\ast=J$ when $k$ is even. In particular we have $J\vert_{V_+^k}=\text{id}_{V_+^k}$, and for any $u\in C^\infty(\Om,V^k_+)$, $\delta J u=\delta u$. Thus, $d+ \delta J$ restricted to $V^k_+$ is elliptic with coerciveness, since it coincides with $d+ \delta$. 
\end{ex}

\begin{rem}
The example of the symplectic Dirac operator can be put in a more general perspective. Assume we are given a homogeneous constant coefficient PDE operator (not necessarily of first order) $\De: C^\infty(\R^n,E)\to C^\infty(\R^n,F)$, that is not elliptic. This implies that $\Ker(\mathbb{D}(\xi))\neq \{0\}$ for some $\xi\in S^{n-1}$. In particular, \emph{the wave cone} of $\De$ satisfies
\begin{align*}
\Lambda_{\De}=\bigcup_{\xi\in S^{n-1}}\Ker(\mathbb{D}(\xi))\neq \{0\}. 
\end{align*}
Note that $\Lambda_{\De}\subset E$ is in general not a linear space. Assume that there exists non-trivial vector space $V\subset (\Lambda_{\De}^c\cup\{0\})$. Then clearly $\mathbb{D}(\xi)\vert_{V}$ is injective and so 
the restricted operator 
\begin{align*}
\De\vert_V: C^\infty(\R^n,V)\to C^\infty(\R^n,F)
\end{align*}
is elliptic. There are however PDE operators so that $\Lambda_{\De}^c\cup\{0\}$ contains no nontrivial linear subspace. This is for example the case with the hyperbolic Dirac operator. Furthermore, since
\begin{align*}
\Lambda_{\De}^c\cup \{0\}&=\Big(\bigcup_{\xi\in S^{n-1}}\Ker(\mathbb{D}(\xi))\Big)^c\cup \{0\}=\bigcap_{\xi \in S^{n-1}}(\Ker(\mathbb{D}(\xi))^c\cup \{0\})
\end{align*}
and $\im(\Aa(\xi)^\ast)\subset \Ker(\mathbb{D}(\xi))^c\cup \{0\}$, it follows that 
\begin{align*}
\bigcap_{\xi \in S^{n-1}}\im(\Aa(\xi)^\ast) \subset \Lambda_{\De}^c\cup \{0\}. 
\end{align*}
Thus, if 
\begin{align*}
V:=\bigcap_{\xi \in S^{n-1}}\im(\Aa(\xi)^\ast)\neq \{0\}, 
\end{align*}
then $\De\vert_{V}$ is elliptic. This implies by \cite[Def. 2.1]{VanS}  that $\De^\ast$ is {\bf not} \emph{canceling}. 
\end{rem}

%========================NEW SECTION==========================================

\section{\sffamily  $\B$-Pseudoconvexity and Generalized Morrey Inequalities}

We begin by proving two different formulas for the Levi 
form \eqref{eq:Levi}.

\begin{Prop}\label{prop:LeviShape}
Let $\Om\subset \R^n$ be a $C^2$-domain and let $\{\kappa_1(x),...,\kappa_{n-1}((x)\}$ be its principle curvatures. Let $\{\mathbf{e}_1(x),...,\mathbf{e}_{n-1}(x)\}$ be a local ON-eigenframe of the shape operator $\mathcal{S}_{\dv \Om}^x$, that is, the operator corresponding to the second fundamental form, of $\dv \Om$. Then the generalized Levi form $\langle u(x),L_B(x)u(x)\rangle$ equals 
\begin{align}
\langle u(x),L_B(x)u(x)\rangle=\sum_{j=1}^{n-1}\kappa_j(x)\vert \mathbb{B}(\mathbf{e}_j(x)) u(x)\vert^2,
\end{align}
for all $u\in C^1(\overline{\Om},F)$ satisfying the boundary condition \eqref{eq:bv2}.
\end{Prop}

\begin{proof}
With notation as in \eqref{eq:nablasplit}, we have
$$
-\mathscr{B}^\ast\mathbb{B}( N(x))=\mathbb{B}(\nu(x))^\ast\dv_\nu \mathbb{B}( N(x))+\mathbb{B}^\ast(\nabla^T)\mathbb{B}( \nu(x)).
$$

Using the boundary condition \eqref{eq:bv2}, the first term gives
$$
  \langle u(x), \mathbb{B}(\nu(x))^\ast \dv_\nu \mathbb{B}( N(x)) u(x)\rangle=
   \langle\mathbb{B}(\nu(x))  u(x), \dv_\nu \mathbb{B}( N(x)) u(x)\rangle=0.
$$
For the second term, we use $\dv_{\mathbf{e}_j(x)}\nu(x)=\kappa_j(x)\mathbf{e}_j(x)$ to get
\begin{align*}
\mathbb{B}(\nabla^T)^\ast\mathbb{B}( \nu(x))&=\sum_{j=1}^{n-1}\mathbb{B}(\mathbf{e}_j(x))^\ast\mathbb{B}( \dv_{\mathbf{e}_j(x)}\nu(x))\\
&=\sum_{j=1}^{n-1}\mathbb{B}(\mathbf{e}_j(x))^\ast\mathbb{B}( \kappa_j(x)\mathbf{e}_j(x))\\
&=\sum_{j=1}^{n-1}\kappa_j(x)\mathbb{B}(\mathbf{e}_j(x))^\ast\mathbb{B}( \mathbf{e}_j(x)).
\end{align*}
Taking the inner product with $u(x)$ now gives the stated formula.
\end{proof}

Next we give a formula for the Levi form in terms of a defining function for $\Om$, which extends the well known formula from several complex variables. 
Recall that $\rho:\R^n \to \R$ is a \emph{defining function} for $\Om$ if
\begin{itemize}
\item[(i)] $\Om=\{x\in \R^n: \rho(x)<0\}$,
\item[(ii)] $(\overline{\Om})^c=\{x\in \R^n: \rho(x)>0\}$, \text{ and}
\item[(iii)] $\nabla \rho(x)\neq 0$ for all $x\in \dv \Om$. 
\end{itemize}
Note that a $C^k$-domain admits a $C^k$-defining function, $k\geq 1$. 
Using a defining function the Gauss map is given by 
\begin{align*}
N(x)=\frac{\nabla \rho(x)}{\vert \nabla \rho(x)\vert}.
\end{align*}

\begin{Lem}   \label{lem:levidefnfcn}
Let $\Om\subset\R^n$ be a bounded $C^2$-domain with a $C^2$ defining function 
$\rho$. 
Then the generalized Levi form $\langle u(x),L_B(x)u(x)\rangle$ equals 
\begin{align*}
\langle u(x),L_B(x)u(x)\rangle=\frac{1}{\vert \nabla \rho(x)\vert}
\sum_{j=1}^{n}\sum_{k=1}^n\frac{\partial^2 \rho}{\partial x_j\partial x_k}
\langle \Bb(e_j) u(x), \Bb(e_k) u(x)\rangle,
\end{align*}
for all $u\in C^1(\overline{\Om},F)$ satisfying the boundary condition \eqref{eq:bv2}.
\end{Lem}

\begin{proof}
A computation gives 
\begin{align*}
\dv_j N(x)=\frac{\nabla \dv_j\rho(x)}{\vert \nabla \rho(x)\vert}-\frac{\nabla \rho(x)}{\vert \nabla \rho(x)\vert^3}\langle \nabla \dv_j\rho(x),\nabla \rho(x)\rangle. 
\end{align*}
The second term will not contribute to 
$\langle u(x),L_B(x)u(x)\rangle$, since $\nabla \rho(x)$ is normal and 
by hypothesis $\Bb(\nu)u(x)=0$.
For the first term, we have
\begin{align*}
\sum_{j=1}^n\Bb(e_j)^\ast \Bb(\dv_j N(x))u(x)&=\frac{1}{\vert \nabla \rho(x)\vert}\sum_{j=1}^n\Bb(e_j)^\ast\Bb(\nabla \dv_j\rho(x))u(x)
\\&=\frac{1}{\vert \nabla \rho(x)\vert}\sum_{j,k=1}^n\Bb(e_j)^\ast\Bb(e_k)\dv_{j,k}^2\rho(x)u(x). 
\end{align*}
Taking the inner product with $u(x)$ now gives the stated formula.
\end{proof}

Extending standard terminology from several complex variables, 
see Example~\ref{ex:Dolbeaut} below,
we make the following definition.

\begin{Def}\label{def:Bsconvex}
Let $\B$ be a first order homogeneous constant coefficent partial differential operator, 
and let $\Om\subset \R^n$ be a bounded $C^2$ domain.  
We say that $\Om$ is {\em strongly} $\B$-{\em pseudoconvex} if
at each $x\in\dv\Om$, the quadratic form
$\langle u(x),L_B(x)u(x)\rangle$ is positive definite on $\ker(\Bb(\nu(x)))$.
\end{Def}

We note that, in contrast to the Laplace coefficients $\mathbb{M}$ in Section~\ref{sec:Weitzen}, the positive definiteness of $\langle u(x),L_B(x)u(x)\rangle$ is not affected if $\B$
is replaced by $B_0 \B$, for some invertible map $B_0\in \text{GL}(E)$. 

Before proving Theorem~\ref{thm:2}, we discuss the relation between 
strict convexity of $\Om$, that is, when all principal curvatures  are positive,
$\inf_{x,j}\kappa_j(x)>0$, and  
strong $\B$-pseudoconvexity of $\Om$, using the following property.

\begin{Def}\label{def:CoCancel}
Let $\B$ be a homogeneous constant coefficent partial differential operator. Define
\begin{align}\label{eq:Cocancel}
K_\B:=\bigcap_{\xi\in S^{n-1}}\Ker(\Bb(\xi))
\end{align}
An operator $\B$  is said to be \emph{cocanceling} if $K_\B=\{0\}$. 
\end{Def}

We recall the observation from \cite[Lem. 3.11]{GR20} that if $\B$ is 
not cocanceling, then restricting $\B$ to act on functions with values
in $K_\B^\perp$, yields an operator which is cocanceling.

\begin{Prop}\label{prop:CoCancel}
If $\B$ is not cocanceling, then there exists no strongly $\B$-pseudoconvex domains.
If $\B$ is cocanceling and non-elliptic, then any $C^2$ strictly convex domain is strongly 
$\B$-pseudoconvex.
\end{Prop}

\begin{rem}
The reason that we require $\B$ to be non-elliptic is that while elliptic operators are always cocanceling, the boundary condition $\Bb(\nu(x))u(x)=0$ implies that $u(x)=0$ for all $x\in \dv \Om$ in the elliptic case. 
\end{rem}
\begin{proof}
Assume that $\B$ is not cocanceling. 
Then there exists $u\ne 0$ such that $\Bb(\xi)u=0$ for all $\xi\in S^{n-1}$.
In particular $\Bb(\nu(x))u= 0$ and $\Bb(\mathbf{e}_j(x))u=0$
for $j=1,\ldots, n-1$ and $x\in \dv \Om$. 
Thus $\langle u(x), L_B(x)u(x)\rangle\equiv 0$, and $L_B(x)$ 
is not positive definite. 

Conversely, assume that $\B$ is cocanceling and that 
$\Om$ is strictly convex.
By Proposition~\ref{prop:LeviShape}, if 
$\Bb(\nu(x))u= 0$ and $\langle u(x), L_B(x)u(x)\rangle= 0$, then
$\Bb(\mathbf{e}_j(x))u=0$ also for $j=1,\ldots, n-1$.
By linearity $\Bb(\xi)u=0$ for all $\xi\in S^{n-1}$, so by 
cocanceling, $u=0$.
\end{proof}

\begin{proof}[Proof of Theorem \ref{thm:2}]
Theorem \ref{thm:Weitzen} shows that
\begin{align}
\int_{\Om}(\vert \A u(x)\vert^2+\vert \B u(x)\vert^2) dx=\int_{\Om}\langle Du(x),\mathbb{M}Du(x)\rangle dx+\int_{\dv \Om}\langle u(x), L_B(x)u(x)\rangle d\sigma(x).
\end{align}
By assumption $\mathbb{M}$ is positive semi-definite and $L_B(x)$ is positive definite on $\dv \Om$, so the generalized Morrey estimate
\begin{align}
\int_{\Om}\vert \A u(x)\vert^2+\vert \B u(x)\vert^2 dx\gtrsim\int_{\dv \Om}\vert u(x)\vert^2 d\sigma(x)
\end{align}
follows.

To deduce the subelliptic estimate \eqref{eq:subE}, we first note the 
inequality
\begin{equation}   \label{eq:Fabesest}
  \|u\|_{W^{1/2, 2}(\Om,F)}^2\lesssim
  \int_\Om |Du(x)|^2 \lambda(x,\dv\Om) dx+ \|u\|_{L^{2}(\Om,F)}^2,
\end{equation}
where $\lambda(x,\dv\Om)$ denotes the distance from $x$ to $\dv\Om$.
This estimate holds without assuming any boundary conditions and for any Lipschitz domain. A proof is in Fabes~\cite{Fabes}. A correction of 
this proof is in the thesis~\cite[Thm. 1.5.10]{AxThesis}. A straightforward 
localization to bounded domains yields \eqref{eq:Fabesest}.

To estimate the square function norm
$\int_\Om |Du(x)|^2 \lambda(x,\dv\Om) dx$, we represent 
$u(x)$ with the generalized Cauchy--Pompeiu integral formula
\begin{equation}   \label{eq:CauchyPomp}
  u(x)=\int_{\dv\Om} E(y-x)\mathbb{D}(\nu(y))u(y)d\sigma(y)- 
  \int_\Om E(y-x) \De u(y) dy= u_1(x)+ u_2(x),
\end{equation}
where $E(x)$ is a fundamental solution, defined as the inverse
Fourier transform of 
$$
  \widehat E(\xi)= -i(\mathbb{D}(\xi)^*\mathbb{D}(\xi))^{-1}\mathbb{D}(\xi)^*.
$$
The proof of this integral formula is completely similar to 
that of \cite[Thm. 8.1.8]{Ros}, but now applying Stokes theorem to
the linear $1$-form
$$
  (y,v)\mapsto E(y-x)\mathbb{D}(v) u(y).
$$
For complete details on the generalized Cauchy--Pompeiu integral formula for general first order operators we refer to the forthcoming work \cite{Duse}.
For the second interior term in \eqref{eq:CauchyPomp}, we estimate
$$
  \int_\Om |Du_2(x)|^2 \lambda(x,\dv\Om) dx
  \lesssim \int_\Om |Du_2(x)|^2 dx
  \lesssim \int_\Om |\De u(x)|^2 dx,
$$
using Plancherel's theorem and noting $|\xi| |\widehat E(\xi)|\lesssim 1$.
For the first boundary term $u_1$ in \eqref{eq:CauchyPomp}, we apply
the $T1$ theorem for square functions as follows.
By applying a partition of unity argument to $u$ on $\dv\Om$, we may 
assume, possibly changing to another ON-basis for $\R^n$, that $u$
is supported on a part of 
$\dv\Om$ which coincides with the graph of a compactly supported $C^2$ function
$x_n=\phi(x_1,\ldots, x_{n-1})$, and that the domain $\widetilde \Om$ above this graph, contains $\Om$.
Changing variables in the first integral in \eqref{eq:CauchyPomp}, we have
$$
  Du_1(x) = \int_{\dv\widetilde\Om}D_x E(y'-x',\phi(y')- \phi(x')-t)\mathbb{D}(\nu(y))u(y',\phi(y'))\sqrt{1+|\nabla\phi(y')|^2} dy',
$$
$x\in\Om$,
writing $y=(y',\phi(y'))$ and $x= (x',\phi(x')+t)$, where $t= x_n-\phi(x')$.
Note that $\lambda(x,\dv\Om)\lesssim t$.
Define the family of kernel functions
$$
  \theta_t(x',y')= t D_x E(y'-x',\phi(y')- \phi(x')-t)\mathbb{D}(\nu(y))
  \sqrt{1+|\nabla\phi(y')|^2}\in \LL(F,F\otimes \R^n).
$$
Let $\theta_t$ be the corresponding integral operator
$(\theta_t u)(x')= \int_{\R^{n-1}} \theta_t(x',y')u(y') dy'$.
It remains to prove the square function estimate
$$
  \int_0^\ell \|\theta_t u\|_{L^2(\R^{n-1},F)}^2 \frac{dt}t
  \lesssim  \| u\|_{L^2(\R^{n-1},F)}^2, 
$$
for $\ell$ comparable to the diameter of $\Om$.
We do this by verifying the hypothesis of the $T1$ theorem for 
square functions. See Christ and Journ\'e~\cite[Thm. 1]{ChrJour}
and its direct vector-valued extension stated by Grau de la Herr\'an
and Hofmann in \cite[Thm. 2.2]{GrauHerrHof}.
It is straightforward to verify the kernel estimates
$$
  |\theta_t(x',y')| \lesssim \frac 1{t^{n-1}}\frac 1{(1+ |x'-y'|/t)^{n}}
$$
and 
$$
  |D_{x'}\theta_t(x',y')|+ |D_{y'}\theta_t(x',y')|
   \lesssim \frac 1{t^{n}}\frac 1{(1+ |x'-y'|/t)^{n}}.
$$
For the latter, recall that we assume that $t\le \ell$ and that 
$\nu$ is constant outside
a compact subset of $\dv\widetilde\Om$.
Finally, we verify the cancellation hypothesis
$$
  \int_{\R^{n-1}} \theta_t(x',y') dy'=0, \qquad\text{for all } x'\in\R^{n-1}, t>0.
$$
Changing variables back to $\dv\widetilde\Om$, this amounts to showing 
\begin{equation}   \label{eq:canc}
  \int_{\dv\widetilde\Om} D_x E(y-x)\mathbb{D}(\nu(y)) d\sigma(y)=0,
\end{equation}
for all $x\in\widetilde\Om$.
To this end, fix a large radius $R$ and define the surfaces
$\Gamma_1= \sett{y\in\dv\widetilde\Om}{|y|<R}$
and $\Gamma_2=\sett{y\in\widetilde\Om}{|y|=R}$.
Applying Stokes' theorem, see \cite[Thm. 7.3.9]{Ros}, 
to the linear $1$-form 
$$
  (y,v)\mapsto E(y-x)\mathbb{D}(v) 
$$
and the domain $\sett{y\in\widetilde\Om}{|y|<R}$, gives
$$
  \int_{\Gamma_1} E(y-x)\mathbb{D}(\nu(y)) d\sigma(y)=1-
  \int_{\Gamma_2} E(y-x)\mathbb{D}(y/R) d\sigma(y).
$$
Applying $D=D_x$, this yields
$$
  \int_{\dv\widetilde\Om, |y|<R} D_xE(y-x)\mathbb{D}(\nu(y)) d\sigma(y)=
  0+ O(1/R),
$$
and \eqref{eq:canc} follows upon letting $R\to \infty$.
This completes the proof of Theorem~\ref{thm:2}.
\end{proof}

While we for the Morrey estimate in Theorem \ref{thm:2} do not a priori require the associated operator $\mathscr{D}$ to be elliptic, to have subelliptic estimates it is necessary that $\mathscr{D}$ is elliptic. Indeed we have.
\begin{Prop}
A necessary condition for a pair $(\A,\B)$ to satisfy a subelliptic estimate on a domain $\Om\subset \R^n$ is that the associated operator $\mathscr{D}$ is elliptic. 
\end{Prop}

\begin{proof}
We argue as in the proof of Theorem \ref{thm:ANS}. Assume that $\mathscr{D}$ is not elliptic. Then there exist $\xi\in \R^n\setminus\{0\}$ and $v\in E\setminus\{0\}$ such that $\Aa(\xi)v=0$ and $\Bb(\xi)v=0$. Let $\eta\in C^\infty_0(\Om)$ and consider 
\begin{align*}
u_t(x)=\eta(x)e^{it\langle x,\xi\rangle}v.
\end{align*}
Then $u_t\vert_{\dv\Om}=0$ and so $u_t$ trivially satisfies the boundary condition $\mathbb{B}(\nu(x))u_t(x)=0$ for $x\in \Om$ and all $t\in \R$. Furthermore, $\mathscr{D}u_t(x)=\mathbb{D}(\nabla \eta(x))e^{it\langle x,\xi\rangle}v$ and $\Vert u_t\Vert_{L^2}\leq C$ and $\Vert \mathscr{D}u_t\Vert_{L^2}\leq C$ for some $C<\infty$ and all $t$. 
But it is clear that
\begin{align*}
\lim_{t\to \infty}\Vert u_t\Vert_{W^{s,2}(\Om)}=\infty,
\end{align*}
for any $s>0$, so  $(\A,\B)$ does not satisfy a subelliptic estimate.
\end{proof}

In \cite[Thm. 1, p 323]{Sw2}, Sweeney gives a necessary and sufficient condition for \eqref{eq:subE} to hold for a short elliptic complex. This relies on first decomposing the complex into a normalized form following a procedure due to Guillemin and Rockland in \cite{Gu} and \cite{Ro} respectively, and then verifying the non-algebraic condition $(G)$ in \cite{Sw2}. However, this condition does not give an explicit geometric characterization of the domains for which subelliptic estimates hold. The sufficient algebraic strong $\B$-pseudoconvexity condition in 
Theorem~\ref{thm:2} generalizes the following classical special case
from several complex variables.

\begin{ex}  \label{ex:Dolbeaut}
   Theorem~\ref{thm:2} generalizes the classical Morrey estimate for the 
   $\dbar$-Neumann problem and the Dolbeaut complex from several complex variables.
   We limit ourselves to $(0,q)$ forms in $\C^n$, $0\le q\le n$,
   and recall that the complex differentials 
   $$d\zbar^\alpha= d\zbar_{\alpha_1}\wedg\ldots \wedg d\zbar_{\alpha_q},$$
   $1\le \alpha_1<\cdots<\alpha_q\le n$,
   form a basis for the space $F=\Lambda^{0,q}$ of $(0,q)$-forms.
   Using the Wirtinger operators 
$$
  \partial_{z_j}= \tfrac 12(\partial_{x_j}-i\partial_{y_j})\quad\text{and}
  \quad \tfrac 12(\partial_{x_j}+i\partial_{y_j}),
$$
in $\C^n= \{(z_j)\}= \{(x_j+iy_j)\}$, 
we consider the operators
$$
  \A f= \dbar f= \sum_{j=1}^n d\zbar_j\wedg \partial_{\zbar_j} f,
$$
and 
$$
  \B f= -\dbar^* f= \sum_{j=1}^n d\zbar_j\lctr \partial_{z_j} f,
$$
acting of $(0,q)$-forms $f$, and where $\A f$ and $\B f$ are $(0,q+1)$
and $(0,q-1)$-forms respectively. 
As noted at the end of \cite[sec. 9.1]{Ros}, $\oplus_q \Lambda^{0,q}$
is isomorphic to the complex spinor space of $\C^n$, with the 
spin--Dirac operator $\slashed{\partial}$ corresponding to the Dolbeaut complex
$\dbar$ in the same way that the Hodge--Dirac operator corresponds to 
the de Rham complex.
As in \cite{Ros}, $\lctr$ denotes the left interior product, that is, the 
pointwise adjoint of left exterior multiplication $\wedg$.

To compute the Laplace coefficients $M_{jk}$, we use real coordinates
$x_j$, $j=-n,\ldots, -1,1,\ldots, n$, for $\C^n= \R^{2n}$, setting $x_{-j}=y_j$, $j=1,\ldots, n$.
Using the anticommutation relation
$$
  d\zbar_j\lctr (d\zbar_k\wedg f)+ d\zbar_k\wedg(d\zbar_j\lctr f)= \delta_{jk}f,
$$
we have $$M_{jk}= \tfrac 14 I_{\Lambda^{0,q}} c_{jk},$$
where
$c_{jk}=0$ when $|j|\ne |k|$, $c_{jj}= c_{-j,-j}=1$
and $c_{j,-j}=i= -c_{-j,j}$, $j=1,\ldots,n$.
It follows that $\mathbb{M}$ is positive semi-definite since 
$\begin{bmatrix} 1 & i \\ -i & 1 \end{bmatrix}$ has eigenvalues 
$0$ and $2$.
Moreover, the operator $\De_M= \sqrt{\mathbb{M}}D= \sqrt 2\mathbb{M}D$ is not $\C$-elliptic. 
Indeed, it is well known that the Dolbeaut complex is not coercive.

To compute the Levi form for the Dolbeaut complex, we sum over $|j|, |k|=1,\ldots, n$ in the real variables expression given in Lemma~\ref{lem:levidefnfcn}. 
Assuming for simplicity that $|\nabla \rho|=1$, and note from the Wirtinger 
operator $\partial_{\zbar_j}$ that 
$$
   \Bb(e_j)u(x)= \tfrac 12 d\zbar_j\lctr u(x)
   \quad\text{and}\quad 
   \Bb(e_{-j})u(x)= -\tfrac i2 d\zbar_j\lctr u(x),
$$
$j=1,\ldots,n$.
Therefore 
\begin{multline*}
  \sum_{|j'|=j}\sum_{|k'|=k}\partial_{x_{j'}}\partial_{x_{k'}}\rho
\langle \Bb(e_{j'}) u(x), \Bb(e_{k'}) u(x)\rangle \\
= \tfrac 14(\partial_{x_j}\partial_{x_k}\rho-i\partial_{y_j}\partial_{x_k}\rho
+i \partial_{x_j}\partial_{y_k}\rho+ \partial_{y_j}\partial_{y_k}\rho)
\langle d\zbar_j\lctr u(x), d\zbar_k \lctr u(x)\rangle \\
=\partial_{z_j}\partial_{\zbar_k}\rho
\langle d\zbar_j\lctr u(x), d\zbar_k\lctr u(x)\rangle,
\end{multline*}
and we arrive at the well known complex variables formula 
$$
\langle u(x),L_B(x)u(x)\rangle=
\sum_{j=1}^{n}\sum_{k=1}^n\frac{\partial^2 \rho}{\partial z_j\partial \zbar_k}
\langle d\zbar_j\lctr u, d\zbar_k\lctr u\rangle=
\sum_{jk\alpha\beta} \epsilon_{jk\alpha\beta}
\frac{\partial^2 \rho}{\partial z_j\partial \zbar_k} u_\alpha \overline{u_{\beta}},
$$
for the classical Levi form, 
where $\epsilon_{jk\alpha\beta}=\pm 1$ is a permutation sign 
when $j\in\alpha$, $k\in\beta$, $\alpha\setminus\{j\}= \beta\setminus\{k\}$,
and else $\epsilon_{jk\alpha\beta}=0$.
See \cite[Secs. 12.2, 12.8]{Taylor}.
Theorem \ref{thm:2} thus recovers the classical result of Morrey that 
we have the subelliptic $1/2$-regularity estimate for the 
Dolbeaut complex under the $\dbar$-Neumann boundary condition
$$
   \sum_{j=1}^n \nu_j(x) d\zbar_j \lctr u(x)=0
$$
on $(0,q)$-forms $u$, for strongly pseudoconvex domains $\Om$.
Moreover, it proves the stronger square function estimate \eqref{eq:subE}.
\end{ex}

\begin{ex}[Dirac type operators]
Assume that $(\A, \B)$ are associated to a short elliptic complex
\begin{align*}   
\xymatrix{C^\infty(\Om,E)  \ar[r]^{\B^*} & C^\infty(\Om,F)\ar[r]^{\A} &C^\infty(\Om,G)},
\end{align*}
meaning that the $\A\circ \B^\ast=0$ and $\text{ker}(\Aa(\xi))=\text{im}(\Bb(\xi)^\ast)$ for all $\xi \in \R^n\setminus\{0\}$. Assume further that the associated operator $\mathscr{D}$ is an elliptic operator of \emph{Dirac type} meaning that 
\begin{align}\label{eq:Dirac1}
\mathbb{D}(\xi)^\ast \mathbb{D}(\eta)+\mathbb{D}(\eta)^\ast \mathbb{D}(\xi)&=\pm2\langle \xi,\eta\rangle \text{id}_{F}\\
\mathbb{D}(\xi) \mathbb{D}(\eta)^\ast+\mathbb{D}(\eta) \mathbb{D}(\xi)^\ast&=\pm2\langle \xi,\eta\rangle \text{id} _{E\oplus G}\label{eq:Dirac2}
\end{align} 
Identities \eqref{eq:Dirac1}-\eqref{eq:Dirac2} expressed in terms of $\Aa$ and $\Bb$ become 
\begin{align}\label{eq:Dirac3}
\Aa(\xi)^\ast \Aa(\eta)+\Aa(\eta)^\ast \Aa(\xi)+\Bb(\xi)^\ast \Bb(\eta)+\Bb(\eta)^\ast \Bb(\xi)&=\pm2\langle \xi,\eta\rangle \text{id}_{F}\\
\Aa(\xi)\Aa(\eta)^\ast+\Aa(\eta) \Aa(\xi)^\ast&=\pm2\langle \xi,\eta\rangle \text{id} _{G}\label{eq:Dirac4}\\
\Bb(\xi) \Bb(\eta)^\ast+\Bb(\eta) \Bb(\xi)^\ast&=\pm2\langle \xi,\eta\rangle \text{id} _{E}\label{eq:Dirac5}\\
\Aa(\xi) \Bb(\eta)^\ast+\Aa(\eta) \Bb(\xi)^\ast&=0\label{eq:Dirac6}\\
\Bb(\xi)\Aa(\eta)^\ast+\Bb(\eta) \Aa(\xi)^\ast&=0
\end{align} 
Let 
\begin{align*}
M(\xi,\eta):=\Aa(\xi)^\ast \Aa(\eta)+\Bb(\eta)^\ast \Bb(\xi).
\end{align*}
Then $M_{ij}=M(e_i,e_j)$, and we see that identity \eqref{eq:Dirac3} can be written as 
\begin{align*}
M(\xi,\eta)+M(\eta,\xi)=2\langle \xi,\eta\rangle \text{id}_{F}
\end{align*}
Since $M(\xi,\eta)^\ast=M(\eta,\xi)$, it follows that $\frac{1}{2}(M(\xi,\eta)+M(\xi,\eta)^\ast)=\pm \langle \xi,\eta\rangle \text{id}_{F}$. 
We now consider the case when $\De=d+\delta$, to which Theorem \ref{thm:1} applies. In this case, with $w\in \Lambda^k\R^n$,
\begin{align*}
\Aa(\xi)w=\xi \wedg w, \quad \Bb(\xi)w=\xi \lctr w,
\end{align*}
and we note that $M(\xi,\eta)=M(\eta,\xi)=\langle \xi,\eta\rangle \text{id}_{\Lambda^k\R^n}$. This property is not true for general operators of Dirac type. Indeed, this property fails for the Dolbeaut complex. 
\end{ex}

\bibliographystyle{alpha}

\end{document}